\documentstyle[12pt]{article}

\newcommand{\be}{\begin{equation}}
\newcommand{\ee}{\end{equation}}
\newcommand{\bn}{\begin{eqnarray}}
\newcommand{\en}{\end{eqnarray}}

\def\one{\raise0.0ex\hbox{$1$\kern-0.35em\raise0.0ex\hbox{$1$}}}
\def\C{\raise0.0ex\hbox{$\bf C$\kern-0.72em\raise0.0ex\hbox{$1$}}}
\def\R{\raise0.0ex\hbox{$\bf R$\kern-0.47em\raise0.0ex
\hbox{$\rule{0.8pt}{1.6ex}$}}}
\def\larrr{\raise0.3ex\hbox{$\longrightarrow$\kern-1.5em\raise-1.1ex
\hbox{$\scriptstyle{r\rightarrow 1}$}}}

\normalbaselines
\begin{document}
\begin{titlepage}

\begin{center}

\vspace*{1.0cm}

{\Large{\bf Geometry of Quantum Group Twists, Multidimensional Jackson
 Calculus and Regularization}}

\vskip 1.5cm

{\large {\bf A.P.Demichev}}\renewcommand{\thefootnote}
{\dagger}\footnote{on leave of absence from
Nuclear Physics Institute,
Moscow State University,
119899, Moscow, Russia}

\vskip 0.5cm

Centro Brasileiro de Pesquisas Fisicas - CBPF/CNPq, \\
Rua Dr.Xavier Sigaund, 150, \\
22290-180, Rio de Janeiro, \\
RJ Brasil

\end{center}

\vspace{3 cm}

\begin{abstract}
\normalsize

We show that R-matricies of all simple quantum groups have the properties which
permit to present quantum group twists as transitions to other coordinate
frames on quantum spaces. This implies physical equivalence of field theories
invariant with respect to q-groups (considered as q-deformed space-time groups
of transformations) connected with each other by the twists. Taking into
account this freedom we study quantum spaces of the special type: with
commuting coordinates but with q-deformed differential calculus and construct
$GL_r(N)$ invariant multidimensional Jackson derivatives. We consider a
particle and field theory on a two-dimensional q-space of this kind and come to
the conclusion that only one (time-like) coordinate proved to be discretized.

\end{abstract}
\end{titlepage}

\section{Introduction}

Lattice regularization has many advantages and plays important role in quantum
field theory (see e.g. \cite{Lattice}). Unfortunately, it has some shortcomings
also. Perhaps the most essential one is a space-time symmetry breaking. The
general reason for the latter is connected with introduction of a lattice in a
theory {\it by hands} and so with its too rigid nature. Among other reasons
this fact initiated many attempts to costruct discrete ("quantized") space-time
manifolds on a more deep background (about old attempts see e.g.
\cite{Blokhintsev}, ch. VII and refs. therein). Last years this problem has
received revival and considerable interest due to appearence of quantum groups
\cite{Drinfel'd}, \cite{Woronowicz}, \cite{ReshetikhinTF}.

Quantum spaces which appear in the frame of quantum group theory
\cite{ReshetikhinTF} have many unusual properties, in particular, q-deformed
differential calculi \cite{WessZ} and, in general, non-commuting coordinates.
In one-dimensional space q-derivative can be represented by Jackson difference
operator \cite{Jackson},\cite{Exton}. This, in turn, provides description of a
quantum patricle on a one-dimensional lattice \cite{DimakisMH}. Thus a
q-deformation of a differential calculus apparently leads to space
discretization. The relation of the non-commutativity of coordinates to space
discretization is not so strightforward. From the other hand, the
non-commutativity causes many problems in construction of field theories on
q-spaces. Indeed, it means that operators of coordinates can not be
diagonalized and have not definite values simultaneously. This essentially new
property of q-spaces is not disastrous for the approach which aims to construct
new and more fundamental theory of space-time, gravity, etc. An interesting
attempts in this direction was made in \cite{Kempf}. From this work it is clear
also that it is special property of phase subspace for each pair of coordinate
and conjugate momentum that leads to improved ultraviolet behaviour of the
model (perturbation theory is based on elementary excitations with fast growing
energy spectrum of q-oscillators, roughly
$
E_n\sim [n;q]\equiv (q^n-1)/(q-1) \ ,
$
where $q>1$). The non-commutativity looks like "payment" for existence of some
q-deformed symmetry ($SU_q(N)$ in this case, $N$ being space-time
dimensionality). As is seen from the paper, the physical properties of the
model are quite unusual. In particular, there are no states with definite
values of position or momentum coordinates and even with definite values of any
part of them. So the very notion of a particle becomes ill defined. This means
that models of this kind can not be considered as a regularization of usual
quantum field theory but prove to be of essentially new type.

These properties of q-deformed models seems to cause serious problems if one
aims to construct in this way lattice-like regularization of {\it existent}
field theories. It is obvious that the necessary condition for the such
regularization (which does not bring essentially new properties to  field
theoretical  models on a continuous space-time) is a commutativity of
coordinates. This letter is devoted to the studing of such possibility. The key
idea of our approach is a physical equivalence of different q-deformations of
space-time manifolds which are related to each other by twists of corresponding
quantum groups. As was shown in \cite{ChaichianD}, q-deformed Minkowski
space-time with non-commuting coordinates which corresponds to pure twisted
Poincar\'e group (i.e. to the q-group obtained from the classical one by a
twist) can be constructed from a usual Minkowski space with help of appropriate
coordinate transformation and q-generalization of 4-beins. In this letter we
generalize this partial result to twists of non-trivially deformed all simple
qroups. More precisely, we will show that known R-matricies for all simple
q-groups have the property which permits to describe the twist procedure as a
transformation of q-space coordinates.

To construct lattice-like regularization one definitely needs multidimensional
finite difference calculus. Using above mentioned freedom in choice of
different q-spaces we consider $GL_q(N)$-invariant q-spaces with commuting
coordinates and q-deformed differential calculus and construct multidimensional
analog of the Jackson calculus (invariant with respect to appropriate quantum
group).

Using explicit formulas for a two-dimensional space we consider a quantum
mechanical particle and simplest field theory and show that the latter is
equivalent to the system on a cylinder with time coordinate  taking values on
equidistant lattice along the cylinder. Surprisingly, the second coordinate
(space-like coordinate on a circle) proved to be continuos. This fact takes its
origin in the properties of involution of the corresponding quantum general
linear group which lead to construction of the quntum mechanical states in a
mixed coordinate-momentum representation, so that "the second discreteness"
corresponds to integer numbers labeling Fourier modes on the circles. Such
discretization is enough for the regularization of two-dimensional theories.
But this fact is dangerous for higher dimensional spaces.

\section{Geometry of quantum group twists}

As is shown in \cite{Drinfel'd3},\cite{Reshetikhin}, multiparametric quantum
groups can be obtained from a one-parametric q-group via so called twists of a
quasitriangular Hopf algebra $A$ with help of an element ${\cal F}=\sum
f^i\otimes f_i \in A\otimes A$, which satisfies certain relations, so that new
coproduct $\Delta^{(F)}$ and new universal R-matrix ${\cal R}^{(F)}$ are
connected with the initial objects $\Delta$ and ${\cal R}$ through the
relations
$$
\Delta ^{(F)}={\cal F}\Delta {\cal F}^{-1}, \qquad {\cal R}^{(F)}={\cal
F}^{-1}{\cal R}{\cal F}^{-1} \ .
$$

Consider at first the case of q-deformations of $GL(N)$ groups. In this case a
twist of R-matrix in fundamental representation $R$ is described with help of
diagonal matrix $F=diag(f_{11},f_{12},...,f_{nn})$ with $f_{ij}f_{ji}=1$ so
that
R-matrix $R^{(F)}$ of the twisted group $GL_{r,\tilde{q}_{ij}}(N)$ has a form
$$
R^{(F)}=F^{-1}RF^{-1} \ .
$$
Here $R$ is (in general, also multiparametric) R-matrix of the initial group
$GL_{r,q_{ij}}(N)$ and
\be
\tilde{q}_{ij}=q_{ij}f^2_{ij}\ .                       \label{1}
\ee
Coordinates of the initial quantum space $C_q^N[x^i]$ satisfy the commutation
relations (CR) \cite{ReshetikhinTF},\cite{Schirrmacher}
\be
x^ix^j=q_{ij}x^jx^i                                     \label{2}
\ee
and coordinates of the twisted space $C^{(F)N}_q[\tilde{x}^i]$ have the CR
\be
\tilde{x}^i\tilde{x}^j=
\tilde{q}_{ij}\tilde{x}^j\tilde{x}^i \ .                    \label{3}
\ee
Now we introduce the algebra $E^N_q[e^i,g_j]$ with the generators
$\{e^i,g_i\}_{i=1}^N$ which commute with coordinates and put
\be
\tilde{x}^i=e^ix^i \qquad\mbox{(no summation)}\ .          \label{4}
\ee
The elements $e^i$ play the role of components of a q-deformed (diagonal)
N-bein. CR for them follows from (\ref{1})-(\ref{4})
\be
e^ie^j=f_{ij}^2e^je^i \ ,                                   \label{5}
\ee
and  $g_i$ are inverse elements
\be
g_ie^i=1\ .                                                      \label{6}
\ee
The coordinates $\tilde{x}^i$ are transformed by a q-matrix $\tilde{T}$
\be
\tilde{x}^{\prime\ i} =\sum_{j=1}^N \tilde{T}^i_{\ j}\tilde{x}^j \ . \label{7}
\ee
Then using (\ref{4}),(\ref{6}) one obtains from (\ref{7}) transformations of
the coordinates $x^i$
\be
x^{\prime\ i}= \sum_{j=1}^N g_i\otimes \tilde{T}^i_{\ j}\otimes e^j\otimes x^j
\ .                                                             \label{8}
\ee
We used in (\ref{8}) a cross-product sign to stress that the elements from the
different sets commute with each other (the elements $g_i$ in (\ref{8}) must be
considered as the inverse elements to generators $e^i$ of another copy of an
algebra $E^N_q$ with respect to the elements $e^i$ entering the same formula).
The relation (\ref{8}) means that the coordinates $x^i$ are transformed by the
matrix $T$ with the entries
\be
T^i_{\ j}= g_i\otimes \tilde{T}^i_{\ j}\otimes e^j \qquad\mbox{(no summation)}
\ .                                                              \label{9}
\ee
Using (\ref{6}) one can express the matrix $\tilde{T}$ through $T$
\be
\tilde{T}^i_{\ j}= e^i\otimes T^i_{\ j}\otimes g_j \qquad\mbox{(no summation)}
\ .                                                              \label{10}
\ee

One can check strightforwardly that $x^{\prime i}$ defined by (\ref{8}) satisfy
the correct CR
$$
x^{\prime i}x^{\prime j}=q_{ij}x^{\prime j}x^{\prime i} \ .
$$

The general reason for this is the following property of the R-matricies: if a
q-matrix $T$ satisfies TT-relation defined by the corresponding R-matrix, then
the $\tilde{T}$-matrix defined by (\ref{9}) or (\ref{10}) satisfies the
relation with twisted R-matrix $R^{(F)}$.

To prove this statement let us write TT-relation (CR for entries of a matrix
$T$) in explicit form
$$
\sum_{p,s}R^{mn}_{\ ps}T^p_{\ u}T^s_{\ v}=\sum_{s,r}T^n_{\ s}T^m_{\ r}
R^{rs}_{\ uv}
$$
and substitute $T^i_{\ j}$ by their expressions (\ref{9}) in terms of
$\tilde{T}^i_{\ j}$. This gives the relation for the latters
\be
\sum_{p,s}R^{mn}_{\ ps}g_pg_se^ue^v\tilde{T}^p_{\ u}\tilde{T}^s_{\ v}=
\sum_{s,r}\tilde{T}^n_{\ s}\tilde{T}^m_{\ r}g_ng_me^se^rR^{rs}_{\ uv}
\ .                                                      \label{11n}
\ee
Remind that in this relation the elements $g_i$ must be considered as inverse
elements for the generators $e^i$ of another copy of an algebra $E^N_q$ and so
they commute with the elements $e^i$ entering the same relation.

Multiparametric R-matrix for $GL_{r,q_{ij}}(N)$ group has a form
\be
R^{mn}_{\ ps}=B^{mn}_{\ ps}+N^{mn}_{\ ps}\ ,           \label{12}
\ee
where $B$ is the diagonal matrix
\be
B^{mn}_{\ ps} = \delta^m_{\ p}\delta^n_{\ s}\bigl(\delta^{mn} +
\Theta^{nm} q^{-1}_{mn} + \Theta^{mn} q_{nm}r^{-1} \bigr) \ , \label{13}
\ee
with $\Theta^{mn} = 1$ if $m>n$, $\Theta^{mn} = 0$ if $m\leq n$, and the matrix
$N$ is off-diagonal part of the R-matrix
\be
N^{mn}_{\ ps}=\delta^m_{\ s}\delta^n_{\ p}\Theta^{mn}\bigl(1-r^{-1}\bigr)
\ .                                                              \label{14}
\ee
Using this expressions one easily obtaines
\be
B^{mn}_{\ ps}g_pg_s=g_mg_nB^{mn}_{\ ps}=
f^{-2}_{mn}g_ng_mB^{mn}_{\ ps}=g_ng_mB^{(F)mn}_{\ ps}\ ,  \label{14a}
\ee
\be
N^{mn}_{\ ps}g_pg_s=g_ng_mN^{mn}_{\ ps} \ ,                 \label{14b}
\ee
so that
$$
R^{mn}_{\ ps}g_pg_se^ue^v=g_ng_me^ue^vR^{(F)mn}_{\ ps}\ ,
$$
where $R^{(F)},B^{(F)}$ are the twisted matricies of the same form
(\ref{12})-(\ref{14}) but for the twisted parameters
$\tilde{q}_{ij}=q_{ij}f^2_{ij}.$ Analogous consideration for the RHS of
(\ref{11n}) shows that this relation can be rewritten in the form
\be
\sum_{p,s}R^{(F)mn}_{\ \ \ \ ps}\tilde{T}^p_{\ u}\tilde{T}^s_{\ v}
=\sum_{s,r}\tilde{T}^n_{\ s}\tilde{T}^m_{\ r}
R^{(F)rs}_{\ \ \ \ uv} \ . \qquad\qquad\Box\
  \label{14c}
\ee

Thus twisted q-matricies can be constructed with the help of q-deformed N-beins
(\ref{5}),(\ref{6}) and the formula (\ref{9}) which is direct generalization
(q-deformation) of a relation between matricies  of transformations in
different coordinate frames.

In the case of q-deformation of simple groups of the series $B_N,\ C_N,\ D_N$
there is one more structure, namely an invariant length \cite{ReshetikhinTF}
$$
L_q=\sum_{i,j}x^iC_{ij}x^j=\sum_{i}l_ix^{i^\prime}x^i\ ,
$$
where $i^\prime=N+1-i$. Values of the coefficients $l_i$ can be found in
\cite{ReshetikhinTF} and are not essential for our consideration. To preserve
$L_q$, components of a q-bein must satisfy the additional constraints
\be
e^ie^{i^\prime}=e^{i^\prime}e^i=1 \ ,            \label{15}
\ee
$i=1,...,N/2$ for $C_N,D_N$ series; $i=1,...,(N+1)/2$ for $B_N$ series. In
particular, for the series $B_N$
\be
e^{(N+1)/2}=1\ .                                  \label{16}
\ee
These constraints reduce  number of twist parameters, which from the
geometrical point of view define CR for the components of the q-beins, so that
the number is equal to $k(k-1)/2$, where $k$ is rank of a group. R-matricies
for the $B_N,\ C_N,\ D_N$ series have the form
$$
R^{ij}_{\ kl} = \biggl[\delta^i_{\ k}\delta^j_{\
l}\biggl(r\delta^{ij}(1-\delta^{ii^\prime}) +
(\Theta^{ji} rq^{-1}_{ij} + \Theta^{ij} q_{ji}r^{-1})(1-\delta^{ii^\prime})
\biggr)+(r-r^{-1})\delta^i_{\ l}\delta^j_{\ k}\Theta^{ij}\biggr]+
$$
$$
+\biggl[\frac{1}{r}\delta^i_{\ k}\delta^j_{\
l}\delta^{ji^\prime}(1-\delta^{ii^\prime})- \Theta^{ij}
(r-r^{-1})r^{(\rho_i-\rho_j)}\epsilon_i\epsilon_j\delta^{ij^\prime}
\delta_{kl^\prime}+\delta^i_{\ (N+1)/2}\delta^j_{\ (N+1)/2}
\delta^{(N+1)/2}_{\ \ \ \ k}\delta^{(N+1)/2}_{\ \ \ \ l}\biggr]
$$
(the last term exists for the $B_N$ series only). Here $\rho_i$ and
$\epsilon_i$ are integer or half integer
numbers \cite{ReshetikhinTF},\cite{Schirrmacher2}. Using this explicit form one
can easily show that the R-matricies have the property analogous to that of the
  $A_N$ groups. Indeed, the terms in the first square brackets have the
structure similar to those of the R-matrix for the $GL_{r,q_{ij}}(N)$ groups.
So they are transformed properly when the elements $e_i,g_j$ move through them
(cf. (\ref{14a}),(\ref{14b})). The terms in the second square brackets are not
changed because of Kronecker symbols $ \delta^{ji^\prime},\ \delta^{ij^\prime}$
or $\delta^{(N+1)/2}_{\ \ \ \ k},\ \delta^i_{\ (N+1)/2}$ and the relations
(\ref{15}),(\ref{16}).

Thus again the matricies $\tilde{T}$ defined by (\ref{9}),(\ref{10}) satisfy
the CR (\ref{14c}) for twisted quantum groups.

\section{q-Spaces with commuting coordinates and multidiminsional Jackson
differential calculus}

It is well known that in one-dimensional case q-deformed differential calculus
can be realized with help of finite difference operation called Jackson
derivative \cite{Jackson},\cite{Exton} which has the form
\be
D_rf(x)=\frac{f(x)-f(rx)}{(1-r)x}
\ ,                                                          \label{17}
\ee
or
\be
\tilde{D}_rf(x)=\frac{f(r^{-1}x)-f(x)}{(r^{-1}-1)x}         \label{17a}
\ee
with the CR
$$
D_rx-rxD_r=1\ ,\qquad \tilde{D}_rx-r^{-1}x\tilde{D}_r=1
$$
(throughout the paper we will assume that $r\leq 1$). In particular, these
derivatives are suitable for description of a quantum mechanichal particle on a
one-dimensional lattice \cite{DimakisMH}. To consider quantum mechanics on
higher dimensional lattices one needs multidimensional generalization of the
Jackson calculus. Such calculus invariant with respect to the q-group
$GL_r(N):=GL_{r,1}(N)$ can be constructed in the space $C^N_r[x^i]$ with
commuting coordinates $x^i$.

As is shown in the previous section, commuting coordinates differ from
non-commuting ones by the non-commuting factors $e^i$. The situation reminds a
transition from usual three-dimensional Euclidian coordinates to well known
quaternions with basis $\{\sigma_i\}_{i=1}^3$: sometimes it is convenient to
put in the coorespondence to the coordinates  $\{x^i\}_{i=1}^3$
non-commutative quaternions  $\hat x^i:=x^i\sigma_i$ (no summation). Though
this transformation brings new algebraic structure and permits to express
three-dimensional rotations in pure algebraic way, the underling geometrical
and physical structure remains the same.

This analogy leads to the conclusion that one can freely choose a most
convenient quantum space among a set of twisted q-spaces. In particular, in the
case of spaces $C^N_{r,q_{ij}}[x^i]$ the most simple choice is the spaces
$C^N_r [x^i]$ with commuting coordinates.The CR for coordinates and derivatives
on this space are the following \cite{Schirrmacher}

\bn
x^ix^j=x^jx^i \ \ \ \forall\ i,j\ ,\quad
\partial_ix^i=1+rx^i\partial_i+(r-1)\sum_{l=i+1}^Nx^l\partial_l\ , \nonumber\\
\partial_i\partial_k =\frac{1}{r}\partial_k\partial_i \ ,\quad
\partial_ix^k=rx^k\partial_i\ ,\quad \partial_kx^i=x^i\partial_k\ ,
\label{18}\\
i<k,\ \ i,k=1,...,N\ . \nonumber
\en

To develop Jackson calculus define the operators of finite dilatations
\be
A_i=1+(r-1)\sum_{j=i}^Nx^j\partial_j                  \label{19}
\ee
with the commutation relations
\be
\begin{array}{cc}
A_iA_k=A_kA_i\ , & \forall\ k,i\ , \\
A_kx^i=x^iA_k\ , & k>i\ , \\
A_ix^k=rx^kA_i\ , & k\geq i\ .
\end{array}                                           \label{20}
\ee
Note that the operators $A_i$ are analogous to operators $Y^i_{\ j}$ of vector
fields on a simple quantum group, introduced in \cite{Zumino}. The  relations
(\ref{19}) permit to express the q-derivatives in terms of $A_i$
\bn
\partial_i=(1-r)^{-1}(x^i)^{-1}(A_{i+1}-A_i)\ , \ \ \ i=1,...,N \label{21} \\
A_{N+1}:=1\ . \nonumber
\en
The relations (\ref{20}) and (\ref{21}) lead to the following realization of
the q-derivatives in a space of functions of $N$ commuting variables
\be
\partial_if(x^1,...,x^N)=\frac{f(x^1,...,x^i,rx^{i+1},...,rx^N)-
f(x^1,...,x^{i-1},rx^i,...,rx^N)}{(1-r)x^i} \ .       \label{22}
\ee
One can check easily that the finite differences (\ref{22}) satisfy the
$GL_r(N)$ invariant CR (\ref{18}) indeed. These differences look like natural
multidimensional generalization of the Jackson derivative (\ref{17}).

As is shown in \cite{WessZ} there are two types of CR for q-derivatives and
coordinates, which are invariant with respect to q-deformed groups. The first
possibility is presented in (\ref{18}), the second one is the following
\bn
\tilde\partial_ix^i=1+\frac{1}{r}x^i\tilde\partial_i+(\frac{1}{r}-1)\sum_{l=1}
^{i-1}x^l\tilde\partial_l\ , \nonumber\\
\tilde\partial_i\tilde\partial_k =\frac{1}{r}\tilde\partial_k\tilde\partial_i \
,\quad
\tilde\partial_ix^k=x^k\tilde\partial_i\ ,\quad
\tilde\partial_kx^i=\frac{1}{r}x^i\tilde\partial_k\ , \label{22a}\\
i<k,\ i,k=1,...,N\ . \nonumber
\en
For these CR it is natural to introduce the operators
\be
\tilde A_i=1+(\frac{1}{r}-1)\sum_{j=1}^ix^j\tilde\partial_j \ , \label{23}
\ee
which commute with each other and have the following CR with the coordinates
\be
\begin{array}{cc}
\tilde A_kx^i=r^{-1}x^i\tilde A_k\ , & k>i\ , \\
\tilde A_ix^k=x^k\tilde A_i\ , & k\geq i\ .
\end{array}                                           \label{24}
\ee
These relations permit to construct the realization of the q-derivatives
$\tilde\partial_i$ in terms of the finite differences
\be
\tilde\partial_if(x^1,...,x^N)=
\frac{f(r^{-1}x^1,...,r^{-1}x^i,x^{i+1},...,x^N)-
f(r^{-1}x^1,...,r^{-1}x^{i-1},rx^i,...,rx^N)}{(r^{-1}-1)x^i} \ ,
\label{25}
\ee
which is the generalization of the Jackson derivative of the second type
(\ref{17a}).

\section{Quantum particle on two-dimensionall quantum space}

In this section we apply the above considered formulas for construction of
quantum mechanics on a two-dimensional quantum plane. For convenience we
rewrite the CR (\ref{18}), (\ref{20}), (\ref{22a}) and  (\ref{24}) in this
particular case denoting: $x^1=z,\ x^2=\bar z,\ \partial_1=\partial,\
\partial_2=\bar\partial,\ \tilde\partial_1=\tilde\partial,\
\tilde\partial_2=\tilde{\bar\partial},\ A_1=A,\ A_2=\bar A,\ \tilde A_1=\tilde
A,\ A_2=\tilde{\bar A}$

\be
z\bar z = \bar z z \ , \quad  \partial\bar\partial =
\frac{1}{r}\bar\partial\partial\ , \quad \tilde\partial\tilde{\bar\partial} =
\frac{1}{r}\tilde{\bar\partial}\tilde\partial\ ,        \label{26}
\ee
\bn
\partial z=1+rz\partial +(r-1)\bar z\bar\partial\ , \qquad
\bar\partial\bar z=1+r\bar z\bar\partial\ , \nonumber \\
\partial\bar z=r\bar z\partial\ , \qquad
\bar\partial z=z\bar\partial\ ,              \label{27}
\en
\be
Az=rzA\ ,\quad A\bar z=r\bar zA\ ,\quad \bar Az=z\bar A\ , \quad
\bar A\bar z=r\bar z\bar A \ ,                      \label{28}
\ee
\bn
\tilde\partial z=1+\frac{1}{r}z\tilde\partial\ , \qquad
\tilde{\bar\partial}\bar z=1+(\frac{1}{r}-1)z\tilde\partial +\frac{1}{r}\bar
z\tilde{\bar\partial}\ , \nonumber \\
\tilde\partial\bar z=\bar z\tilde\partial\ , \qquad
\tilde{\bar\partial} z=\frac{1}{r}z\tilde{\bar\partial}
\ ,                                                    \label{29}
\en
\be
\tilde Az=\frac{1}{r}z\tilde A\ ,\quad \tilde A\bar z=\bar z\tilde A\ ,\quad
\tilde{\bar A}z=\frac{1}{r}z\tilde{\bar A}\ , \quad
\tilde{\bar A}\bar z=\frac{1}{r}\bar z\tilde{\bar A}\ ,    \label{30}
\ee
All A-operators commute with each other and satisfy the relations
\be
\tilde{\bar A}A=1\ ,\qquad \tilde{\bar A}\bar A=\tilde A\ ,\qquad \tilde
AA=\bar A\ .                                                    \label{31}
\ee
The simplest way to find them is to derive from an action of the operators on
an arbitrary function $f(z,\bar z)$.
The CR between different types of the q-derivatives have the form
\be
\begin{array}{ccc}
\tilde\partial\partial = r\partial\tilde\partial\ , & \qquad &
\tilde{\bar\partial}\bar\partial = r\bar\partial\tilde{\bar\partial}\ ,\\
\tilde\partial\bar\partial = \bar\partial\tilde\partial\ , & \qquad &
\tilde{\bar\partial}\partial = r^2\bar\partial\tilde{\bar\partial}\ .
\end{array}                                                    \label{31a}
\ee

Now we must define $*$-involution in the algebra of the operators which enter
the relations (\ref{26})-(\ref{31a}).

We want to consider the parameter $r$ as a lattice spacing, hence, it must be a
real number. The appropriate involution for $GL_r(2)$ in this case is the
following \cite{Schirrmacher2}
$$
T^*=CTC\ , \qquad
C=\left(
\begin{array}{cc}
0 & 1 \\
1 & 0
\end{array}\right) \ .
$$
This means that
$$
T^*=\left(
\begin{array}{cc}
a^* & b^* \\
c^* & d^*
\end{array} \right) =
\left(
\begin{array}{cc}
d & c \\
b & a
\end{array} \right)
$$
and $z^*=\bar z$. It is not difficult to see that the CR for the A-operators
and the coordinates are consistent with the involution
\be
A^*=\tilde{\bar A}\ , \qquad \bar A^*=\tilde A   \ .         \label{32}
\ee
This gives the involution rules for the derivatives
$$
\partial^*=-\tilde{\bar\partial}+(\frac{1}{r}-1)\frac{z}{\bar z}\tilde\partial
+\frac{1}{\bar z}\ , \qquad
\bar\partial^*=-\tilde\partial+\frac{1}{z}\ .
$$
These rules look rather combersome, but they are direct generalization of the
involution in one-dimensional case \cite{DimakisMH}.

The next step is to construct the representation of the operators in a Hilbert
space so that the involution would coincide with hermitian conjugation
(therefore in the following we will denote the involution by a sign of
hermitian conjugation). To construct convenient basis of a Hilbert space one
needs hermitian operators. Hermitian combinations of the coordinates of the
form
$$
x=(z+\bar z)/2\ , \qquad y=(z-\bar z)/2i
$$
are not convenient as they have identical CR with a part of A-operators and
rather combersome CR with another part. A-operators play the role of conjugate
momenta on a lattice (cf. \cite{DimakisMH}). A natural choice of position
operators follows from the observation that the A-operators generate finite
dilatations rather than translations. This implies the introduction of the
hermitian operator
\be
\rho=\sqrt{\bar z z}\ , \qquad \rho^\dagger=\rho\ ,        \label{33}
\ee
and the unitary one
\be
\Phi=\sqrt{\bar z z^{-1}}\ , \qquad \Phi^\dagger\Phi=1\ .   \label{34}
\ee
As follows from (\ref{32}) and (\ref{31}) we have also unitary operator $A$ and
hermitian operator $B:=\tilde A\bar A$. These operators have the following CR
$$
\begin{array}{ccc}
A\rho = r\rho A\ , & \qquad &
B\rho = \rho B\ ,\\
A\Phi = \Phi A\ , & \qquad &
B\Phi =r\Phi B \ .
\end{array}
$$
All other operators can be expressed in terms of these four ones. Thus we have
two commuting hermitian operators $\rho,\ B$ and can define a basis $\mid N,m>$
in a Hilbert space ${\cal H}$ in which these operators are diagonal (cf.
\cite{DimakisMH}, \cite{KempfM})
\be
\begin{array}{ccc}
\rho\mid N,m>  & = &  \rho_0r^{-N}\mid N,m> \\
B\mid N,m> & = & b_0r^{-m}\mid N,m>\ .
\end{array}                                          \label{35}
\ee
The constants $\rho_0$ and $b_0$ mark different equivalent representations and
for shortness we will put $\rho_0=b_0=1$.
The unitary operators transform a vector with given eigenvalue into the  one
with a neighboring eigenvalue
\be
\begin{array}{ccc}
A\mid N,m>  & = &  \mid N+1,m> \\
\Phi\mid N,m> & = & \mid N,m-1>\ .
\end{array}                                          \label{36}
\ee
So starting from $GL_r(2)$ invariant differential calculus on a quantum plane
we have come naturally to the polar coordinates, the operator $\rho$ being the
operator of radial coordinates and the operator $\Phi$ being of the form $\Phi
= e^{-i\phi}$, where $\phi$ is (hermitian) operator of an angle coordinate. The
structure of the algebra involution leads to the mixed representation with one
coordinate ($\rho$) and one momentum diagonal operators. In classical case a
multiplication of a function by $\Phi$ shifts its Fourier components numbers by
minus unity. This corresponds to the action (\ref{36}) of the operator $\Phi$
on vectors of ${\cal H}$. This implies that $B$ is connected with an angular
momentum operator. In fact it is analog of the operator of the form
\be
\exp\{i\phi_0\frac{\partial}{\partial\phi}\}=\exp\{\phi_0 M\}   \label{36a}
\ee
in the case of quantum mechanics on usual continuous plane, where $\phi_0$ is
some fixed angle and $i\partial/\partial\phi$ is (two-dimensional) angular
momentum operator. Eigenfunctions of this operator are periodic functions
$\exp\{in\phi\}$ with eigenvalues $\exp\{n\phi_0\}=(e^{\phi_0})^n$. The latter
expression coincides with eigenvalue of operator $B$ if one equates
$r=\exp\{\phi_0\}$.

Thus the operator $\rho$ defines values of radial coordinate and operator the
$A$ shifts them (play the role of conjugate momentum). Analogously, the
operator $B$ defines angular momentum values and the operator $\Phi$ of the
conjugate coordinate shifts its eigenvalues.

Now we can consider a q-subgroup $\Lambda$ of the $GL_r(2,\R\ )$ of matricies
of the form
$$
T=\left(
\begin{array}{cc}
a & 0 \\
0 & a^*
\end{array}\right) \ ,
$$
with $aa^*=a^*a$. Because of the latter relation it looks like ordinary group
isomorphic to a multiplicative group of complex numbers. But one must remember
about CR with the generators $Y^i_{\ j}$ of the corresponding quantum universal
enveloping algebra. For left-invariant and right-covariant generators they have
the following general form \cite{Zumino}
\be
Y^i_{\ j}T^k_{\ s}=T^k_{\ l}Y^m_{\ n}(\tilde{R}_{21})^{il}_{\ mt}
(\tilde{R}_{12})^{nt}_{\ js}\ ,                           \label{37}
\ee
where $\tilde{R}_{12}$ and $\tilde{R}_{21}$ are properly normalized R-matricies
of the $GL_r(2,\R\ )$. The CR (\ref{37}) gives for the subgroup $\Lambda$ (i.e.
if $b=c=0$)
\be
\begin{array}{ccc}
Y^1_{\ 1}a = raY^1_{\ 1}\ , & \qquad &
Y^1_{\ 1}a^* = a^* Y^1_{\ 1}\ ,\\
Y^2_{\ 2}a = a Y^2_{\ 2}\ , & \qquad &
Y^2_{\ 2}a^* =ra^* Y^2_{\ 2} \ .
\end{array}                                          \label{38}
\ee
It is easy to see that CR for ${\cal D}:=Y^1_{\ 1}Y^2_{\ 2},\
\bar{\cal D}:=Y^2_{\ 2}, \ a$ and $a^*$ are the same as for the $A,\bar
A,z,\bar z$. So the algebra (\ref{28}) on the quantum plane is isomorphic to
that of the q-subgroup $\Lambda$ and we can identify the q-plane with this
subgroup. The operators $\tilde{A},\tilde{\bar A}$ correspond to
right-invariant and left-covariant generators (cf. \cite{Zumino}). From the
other hand, q-subgroup $\Lambda$ can play the role of symmetry group, the whole
$GL_r(2,\R\ )$ group being the group of linear canonical transformations (of CR
(\ref{27}),(\ref{29})).

The coaction of $\Lambda$ on the coordinates is
\be
\begin{array}{ccc}
z\rightarrow z^\prime = a\otimes z\ , & \quad &
\bar z\rightarrow {\bar z}^\prime = a^*\otimes \bar z\ , \\
\rho\rightarrow \rho^\prime =a_{\rho}\otimes \rho\ , & \quad &
\Phi\rightarrow \Phi^\prime =a_{\Phi}\otimes \Phi\ , \\
a_{\rho}=\sqrt{a^*a}\ , & \quad & a_{\Phi}=\sqrt{a/a^*} \ .
\end{array}                                         \label{39}
\ee
The coordinates $\rho^\prime,\ \Phi^\prime$ have a representation in a Hilbert
space ${\cal H\otimes H}$ (cross product of the same Hilbert space as in our
case the comodule coincides with the symmetry q-group) which has a subspace
(diagonal) isomorphic to ${\cal H}$. The latter corresponds to the
representation of CR based on the coordinates $\rho^\prime,\ \Phi^\prime$ as a
primary ones. All A-operators and so the operator $B$ are invariant with
respect to the coaction (\ref{39}).

To consider a two dimensional quantum free partical one can use, for example,
the Hamiltonian
\be
h=\frac{\Omega}{(1-r)^2}\left[ (1-A)(1-A^\dagger)+(1-B)^2\right]\ ,\label{39aa}
\ee
where $\Omega$ is constant with dimension of energy. In order to find
eigenfunctions and eigenvalues of the Hamiltonian introduce the operators
$$
\rho^{iP}=\exp\{iP\ln\rho\}\ .
$$
and the states (cf. \cite{KempfM})
$$
\mid\one ,m>=\sum_{m=-\infty}^{\infty}\mid N,m>\ ,
$$
$$
\mid P,m>=\rho^{iP}\mid\one,m>\
$$
with the properties
$$
A\mid\one ,m>=\mid\one ,m>\ ,
$$
$$
A\mid P,m>=r^{iP}\mid P,m>\ ,
$$
where $P$ is a real number: $0\leq P\leq\pi/\chi$,\  $\chi := \ln r$. For
eigenvalues of the Hamiltonian (\ref{39aa}) we obtain
$$
h\mid P,m>=\xi_{P,m}\mid P,m>\ ,
$$
$$
\xi_{P,m}=\Omega\left(\frac{(1-\cos\chi P)}{(1-r)^2}+[m;r]^2\right)
\larrr
$$
$$
\larrr \ \Omega (P^2+m^2)\ , \qquad [m;r]=\frac{1-r^m}{1-r} \ .
$$
Thus the operator $h$ has the correct continuous limit but its eigenvalues are
not invariant with respect to the reflection $m\rightarrow -m$. This means that
left and right modes have different properties and positive modes have
decreasing spectrum which can lead to additional divergencies in corresponding
field theories. These properties are caused by exponent-like form of the
operator $B$ analogous to (\ref{36a}) as we discussed above. So it is natural
to consider the Hamiltonian
$$
H=\frac{\Omega}{(1-r)^2}\left[ (1-A)(1-A^\dagger )+\ln^2 B\right]
$$
with the eigenvalues $\Omega\lambda_{P,m}$, where
\be
\lambda_{P,m}=\left(\frac{(1-\cos\chi
P)}{(1-r)^2}+\frac{\chi^2}{(1-r)^2}m^2\right)\larrr \ (P^2+m^2)\ .
\label{39a}
\ee

To construct two dimensional quantum field theory we need a kind of integral
over the variable $\rho$. It can be defined analogously to that in
\cite{DimakisMH}, \cite{KempfM}.

For m-th component $f_m(\rho )$ of "Fourier expansion" we define
\bn
\int_0^Kd_r\rho\ f_m(\rho ):=<K,m\mid\partial^{-1}_{\rho}f_m(\rho )\mid\one ,m>
\nonumber \\
=(1-r)<K,m\mid (1-A)^{-1}\rho f_m(\rho )\mid\one ,m>
\nonumber \\
=(1-r)\sum_{l=-\infty}^{K}r^{-l}f_m(r^{-l})\ ,           \label{40}
\en
where the derivative $\partial_\rho$ is defined as in (\ref{21})
$$
\partial_{\rho}:=\frac{1}{(1-r)\rho}(1-A)\ .
$$
The last expression for the integral in (\ref{40}) has the form of a usual
Jackson integral \cite{Exton}.

To construct the action, note that the two-dimensional space with the
coordinates $z,\bar z$ and the symmetry transformations (\ref{39}) can be
considered as a result of the conformal map
\be
u\rightarrow z=e^u\ , \qquad \bar u\rightarrow \bar z=e^{\bar u}  \label{40a}
\ee
of a cylinder with the coordinates $\bar u,u$. This map is well known and is
used widely, in particular, in conformal and string theories (in frame of so
called radial quantization, see e.g.\cite{GreenSW}).
A coordinate along the cylinder is associated with a time coordinate and space
coordinate takes values on a circle. After the conformal mapping the coordinate
$\rho$ plays the role of time and $\Phi$ of space coordinate. It is easy to see
that in coordinates $u,\bar u$ the transformations (\ref{39}) becomes
translations and the time coordinate $\tau :=Re\ u$ takes values on an
equidistant lattice with a spacing $\ln r$. In continuous case a field theory
free action for a scalar field $\Psi$ in polar coordinates on the $z$-plane has
the form
\be
S_0=\int_o^\infty \frac{d\rho}{\rho}d\phi\ \Psi (\rho ,\phi )
(\rho\partial_{\rho}\rho\partial_{\rho}+\partial^2_{\Phi})\Psi(\rho ,\phi )\ .
                                                         \label{41}
\ee
This expression is explicitly invariant with respect to dilatations of $\rho$
and translations of $\phi$ (continuous analog of (\ref{39})).

A quantum plane analog has the following form
$$
S^r_0=\lim_{K\rightarrow\infty}\sum_{m=-\infty}^{\infty}<K,m\mid
(1-A)^{-1}\Psi_{-m}(\rho )G^{-1}\Psi_m (\rho )\mid \one ,m>\ .
$$
Here the operator $G^{-1}$ has the form similar to the Hamiltonian $H$
$$
G^{-1}=\frac{H}{\Omega} \ ,
$$
though its meaning is quite different, of course.
The finite difference operator $G^{-1}$ in this action has the eigenvalues
$\lambda_{P,m}$ presented in (\ref{39a}). In the case of interacting $\Psi^4$
theory a one-loop correction to, for example, mass term is proportional to a
trace of the operator $G$
$$
\sum_{m=-\infty}^\infty \int_0^{-\pi /\ln r}\frac{dP}{\lambda_{P,m}} \ ,
$$
which has no divergencies. Note, however, that summation over orbital number
$m$ is the same as in continuous case.

\section{Conclusion}

The main result of this work are the relations (\ref{4}) and (\ref{9}) which
show that twists of quantum groups and corresponding q-spaces can be realized
with help of the auxiliary non-commuting elements $e^i,g_k$ satisfying the
relations (\ref{5}) and (\ref{6}). Geometrically this can be interpreted as a
transition to another coordinate frame on a q-space and so for many problems
twisted quantum groups must be physically equivalent. Let us mention again the
analogy between the algebra of $e^i,g_k$ and well known quaternions $\sigma_i$:
one can use either usual three-dimensional coordinates $x^i$ or non-commutative
quaternions $\hat x^i:=x^i\sigma_i$ (no summation). Physical content of a
problem is not changed. Unfortunately, many important quantum groups, for
example, orthogonal ones, have not twisted counterparts with commuting
coordinates of a corresponding quantum space. Note, however, that linear groups
play important role as groups of space-time symmetries. It is enough to mention
$SL(2,C)$ as universal covering of Lorentz group and $SU(2,2)$ as covering of
four-dimensional conformal group.

Using q-deformed space with commuting coordinates we constructed
multidimensional q-group invariant generalization of the famous Jackson
derivatives.

As another application of q-spaces with commuting coordinates we considered
quantum mechanics and simple field theory on a two-dimensional quantum space.
The structure of the involution leads to the mixed coordinate-angular momentum
representation of states of the system. This, in turn, results in
discretization of only one (radial) coordinate of the space in spite of the
q-deformed differential calculus, the spectrum of the angular momentum operator
being unbounded as in usual continuous case. Such partial discretization of a
space-time is enough for two-dimensional models but can be dangerous for higher
dimensional cases. Another lesson from the considered models is that
q-derivatives or finite differences constructed with help of operators of the
form (\ref{19}) are connected with dilatations and not with translations. This
implies that corresponding coordinates are related with usual ones (in which
derivatives generate translations) by non-linear exponent-like map of the type
(\ref{40a}). It is clear, that symmetry transformations, e.g. of Lorentz group,
have quite different form in non-linearly transformed coordinate frames. This
remark can be important for numerous attempts to construct quantum deformations
of relativistic symmetry.

\end{document}